\documentclass[12pt,a4paper]{article}
\usepackage{amsfonts}
\usepackage{amsmath}
\usepackage{mathrsfs}
\usepackage{graphicx}
\usepackage{amsmath,amsthm,amssymb,amscd}
\newtheorem{Thm}{Theorem}[section]
\newtheorem{Lem}[Thm]{Lemma}
\newtheorem{Prop}[Thm]{Proposition}
\newtheorem{Def}[Thm] {Definition}

\newtheorem{Cor}[Thm]{Corollary}
\theoremstyle{remark}
\newtheorem{Rem} [Thm]{Remark}

\theoremstyle{claim}

\newtheorem{Que}[Thm]{Question}

\begin{document}

\begin{center} {\large \bf Non-uniform Hyperbolicity and Non-uniform Specification}
\end{center}
\smallskip

\begin{center}
  Krerley Oliveira$^{*}$\\

  Instituto de Matem$\acute{\textrm{a}}$tica, UFAL 57072-090 Macei$\acute{\textrm{o}}$, AL, Brazil \\

E-mail: {\it krerley@gmail.com}\\
\bigskip

  Xueting Tian$^{\dagger}$\\

  Academy of Mathematics and Systems Science, Chinese Academy of Sciences, Beijing 100190, People's Republic of  China\\ $\&$ \\
Instituto de Matem$\acute{\textrm{a}}$tica, UFAL 57072-090 Macei$\acute{\textrm{o}}$, AL, Brazil \\

  E-mail: {\it tianxt@amss.ac.cn $\,\,\,\&\,\,\,$ txt@pku.edu.cn} \\
\end{center}

\footnotetext {$^{*}$ Krerley is supported by CNPq, CAPES, FAPEAL, INCTMAT and PRONEX.}
\footnotetext
{$^{\dagger}$ Tian is the corresponding  author and supported by CAPES.}
\footnotetext{AMS Review: 37A35;  37D05; 37C35 }

\def\abstractname{\textbf{Abstract}}

\begin{abstract}\addcontentsline{toc}{section}{\bf{English Abstract}}
In this paper we deal with an invariant  ergodic hyperbolic measure $\mu$ for a  diffeomorphism $f,$ assuming that $f$ it is either $C^{1+\alpha}$ or $f$ is $C^1$ and the Oseledec splitting of $\mu$ is dominated. We show that this system $(f,\mu)$  satisfies a weaker and non-uniform version of specification,  related with notions studied in several recent papers, including \cite{STV,Y, PS, T,Var, Oli}.

Our main results have several consequences: as corollaries, we are able to improve the results about quantitative Poincar\'e recurrence, removing the assumption of the non-uniform specification property in the main Theorem of \cite{STV} that establishes an inequality between Lyapunov exponents and local recurrence properties. Another consequence is the fact that any of such measure is the weak limit of averages of Dirac measures at periodic points, as in \cite{Sigmund}. Following \cite{Y} and \cite{PS}, one can show  that the topological pressure can be calculated by considering the convenient weighted sums on periodic points, whenever the dynamics is positive expansive and every measure with pressure close to the topological pressure is hyperbolic.
\end{abstract}

\section{Introduction} \setlength{\parindent}{2em}

In seminal works, Bowen \cite{Bow,Bowen2} and Sigmund \cite{Sigmund}  introduced in the 1970's the notion of \emph{specification} and used this to show many ergodic properties for dynamical systems satisfying this property, including subshifts of finite type and sofic subshifts,  the restriction of an axiom A diffeomorphism to its non wandering set, expanding differentiable maps, and  geodesic flows on manifolds with negative curvature. Before continue, let us recall the definition of the \emph{uniform} specification property (SP) as defined by Bowen in \cite{Bow}.

\newpage

We say that $f: X\rightarrow X$ has the \emph{specification property} (SP) if given $\theta>0$ there exists $K_\theta \in \mathbb{N}$ such that for all $x \in X$, there exists $p \in \mathbb{N}$ such that the dynamical ball

$$B^n_m(x,\theta):=\bigcap\limits_{k=-m}^nf^{-k}\big(B(f^k(x),\theta)\big)$$
contains a periodic point with period $p\leq n+m+K_\theta.$

While the specification property usually holds in symbolic dynamics, often in a trivial way, it is a strong hypothesis for dynamical systems. A very interesting result shows that for continuous $f: [0,1]\rightarrow [0,1]$,  specification holds if and only if $f$ is topologically
mixing (see \cite{Bl83}). On the other hand, it is known that for $\beta$-shifts, specification occurs rarely: it is verified only on a set
of $\beta$ of Lebesgue measure zero (see \cite{Bu97}).

For maps that do not have specification, several weaker versions  were introduced in the literature along the last decade. From the best of our knowledge, the first weaker version was introduced by Marcus at \cite{Mar}, to show that periodic measures are weakly dense for  every ergodic toral automorphism, extending a previous work of Sigmund that proved the same for hyperbolic toral automorphisms. For the non-uniformly hyperbolic case, several versions of non-uniform specification, including \cite{Hir,LLS,LST,ST},   were introduced and  had been used to show many ergodic properties.

\medskip

Saussol,  Troubetzkoy and Vaienti in \cite{STV} also introduced another non-uniform specification property and showed that every hyperbolic ergodic measure having non-uniform specification property satisfies an inequality between Lyapunov exponents and local recurrence rates. This notion is defined using \emph{slow varying} functions that appear in the Pesin's Theory, i.e., functions $q$ that satisfy $q(f^{\pm1}(x))\leq e^\eta q(x)$ for every $x \in M$. Let us define the non-uniform specification precisely.

\begin{Def}\label{Def:NS}

Let $\mu$ be an invariant measure of $f$, we say that $(f,\mu)$ satisfies non-uniform specification property (NS), if for $\mu$ almost every $x$, any small $\eta>0,$ any $\eta-$slowing varying positive function $q$,   any integers $m,\,n,$ and any $\theta>0$ there exists $K:=K(\eta,\theta,x,m,n)$ such that:\\
(i) the non-uniform dynamical ball $$\tilde{B}^n_m(x,\theta):=\cap_{k=-m}^nf^{-k}B(f^k(x),\theta q(f^k(x))^{-2})$$
contains a periodic point with period $p\leq n+m+K;$
\\
(ii) the dependence of $K$ on $m,\,n$ satisfies
$$\lim_{\eta\rightarrow 0}\limsup_{m,n\rightarrow+\infty}\frac{K(\eta,\theta,x,m,n)}{m+n}=0.$$

In words: pieces of orbits with length $n$ are shadowed by periodic orbits whose period is less than $n$ plus a sub linear term with respect to $n$.

\end{Def}
\bigskip

\newpage

  A natural question arises: whether this non-uniform specification property is  valid for non-uniformly hyperbolic systems?

\medskip

In \cite{Oli}, the first author  had proved that every expanding ergodic measure for maps with non-flat critical set, i.e., a strongly mixing measure with positive exponents,   has this non-uniform specification property. In this paper, we give a positive answer for hyperbolic ergodic measures to show that every hyperbolic ergodic measure naturally has the non-uniform specification property of \cite{STV} and thus we can remove the assumption of non-uniform specification property in the main theorem of \cite{STV} to establish an inequality between Lyapunov exponents and local recurrence properties and some other corollaries. Now we state our main results.

\begin{Thm}\label{Thm:NS} Let $f$ be a $C^{1+\alpha}\,(\alpha>0)$ diffeomorphism. Then, every $f$ ergodic hyperbolic measure $\mu$ satisfies NS.
\end{Thm}

\begin{Rem}\label{Rem-NS} In particular, if $\mu$ is a mixing hyperbolic measure, then we have a stronger version of the non-uniform specification property: for any $p\geq n+m+K,$ the non-uniform dynamical ball $\tilde{B}^n_m(x,\theta)$ contains a periodic point with period $p$.
\end{Rem}

\bigskip

\begin{Rem}\label{Rem-BowenBall} Note that if the function $q(\cdot)\equiv 1$, then the non-uniform dynamical ball $\widetilde{B}_m^n(x,\theta,q)$ is the general well-known dynamical ball $B_m^n(x,\theta)$. So, in the conclusion of  non-uniform specification property we can replace non-uniform dynamical balls by dynamical balls and $\eta$ can be omitted.
\end{Rem}

%
%
%
%
%
%
%\begin{Thm}\label{Thm:NS-BowenBall}
% Let $f$ be a $C^{1+\alpha}\,(\alpha>0)$ diffeomorphism. Then every $f$ ergodic hyperbolic measure $\mu$ satisfies non-uniform specification property for dynamical balls in the following sense. For $\mu$ almost every $x$, any integer $m,\,n,$ and any $\theta>0$ there exists $K:=K(\theta,x,m,n)$ such that:\\
%(i) the non-uniform dynamical ball $$B^n_m(x,\theta):=\cap_{k=-m}^nf^{-k}B(f^k(x),\theta)$$
%contains a periodic point with period $p\leq n+m+K$(if $\mu$ is mixing, this property can be hold for any $p\geq n+m+K$);
%\\
%(ii) the dependence of $K$ on $m,\,n$ satisfies
%$$\limsup_{m,n\rightarrow+\infty}\frac{K(\theta,x,m,n)}{m+n}=0.$$
%\end{Thm}
%

\medskip
Given $x\in M$ and $r>0$, denote the first return time of a ball $B(x,r)$ radius $r$ at $x$ by $$\tau(B(x,r)):= \min \{k>0 \,|\,\,f^k(B(x,r))\cap B(x,r)\neq \emptyset\}.$$ Then we can obtain a corollary as follows by using Theorem \ref{Thm:NS} and the main Theorem in \cite{STV}.

\begin{Cor}\label{Cor:Recurrece-Lyaexp}

Let $f$ be a $C^{1+\alpha}\,(\alpha>0)$ diffeomorphism. Then for every $f$ ergodic hyperbolic measure $\mu$,  one has for $\mu$ a.e. $x\in M$, $$\limsup_{r\rightarrow 0}\frac{\tau(B(x,r))}{-\log r}\leq \frac1{\lambda_u}-\frac1{\lambda_s},$$ where $\lambda_u,\,\lambda_s$ are the minimal positive Lyapunov exponent and maximal negative Lyapunov exponent of $\mu$, respectively.
\end{Cor}

\medskip

\begin{Rem}\label{Rem-Rec-Lyexp}  From the main Theorem in \cite{STV}, if $h_{\mu}(f)>0$, then  for $\mu$ a.e. $x\in M$, $$\liminf_{r\rightarrow 0}\frac{\tau(B(x,r))}{-\log r}\geq \frac1{\Lambda_u}-\frac1{\Lambda_s},$$ where $\Lambda_u,\,\Lambda_s$ are the maximal positive Lyapunov exponent and minimal negative Lyapunov exponent of $\mu$, respectively. Therefore,  using our Corollary \ref{Cor:Recurrece-Lyaexp} we can get that if $M$ is two dimensional and $h_{\mu}(f)>0$, then for $\mu$ a.e. $x\in M$,
$$\lim_{r\rightarrow 0}\frac{\tau(B(x,r))}{-\log r}= \frac1{\lambda_u}-\frac1{\lambda_s}=\frac1{\Lambda_u}-\frac1{\Lambda_s}.$$
\end{Rem}

\bigskip

We observe that the specification property that we just discussed before has many similar different forms, as it is presented by \cite{Y}. Let us discuss a slightly generalization of the non-uniform specification above, that we call generalized non-uniform specification property with respect to
several orbit segments(being a generalization of NS introduced in Definition \ref{Def:NS}).

\begin{Def}(Generalized Non-uniform Specification GNS) We say that $\mu$ has the \emph{generalized non-uniform specification property} if  for $\mu$ almost every $x$, any small $\eta>0,$ any $\eta-$slowing varying function $q$ (that is, $q(f^{\pm1}(x))\leq e^\eta q(x)$),   any integer $m,\,n,$ and any $\theta>0$ there exists $K:=K(\eta,\theta,x,m,n)$ satisfying
$$\lim_{\eta\rightarrow 0}\limsup_{m,n\rightarrow+\infty}\frac{K(\eta,\theta,x,m,n)}{m+n}=0$$ and so that the following holds: given points $x_1,x_2,\cdots,x_k$ in a full $\mu$-measure set and positive integers $m_1,\cdots,m_k,n_1,\cdots,n_k,$ there is
$ p_i\geq 0 $ with $$\sum_{i=1}^{k} p_i\leq \sum_{i=1}^{k} K(\eta,\theta,x_i,m_i,n_i)$$
$\big{(}$in particular if $\mu$ is mixing, for any  $t\geq \sum_{i=1}^{k} K(\eta,\theta,x_i,m_i,n_i)$ there is $p_i$
with $\sum_{i=1}^{k} p_i=t$$\big{)}$
%$0\leq p_i\leq K(\eta,\theta,x_i,m_i,n_i)$(in particular if $\mu$ is mixing, $p_i$
%can be chosen arbitrary integer $\geq K(\eta,\theta,x_i,m_i,n_i)$)
and a periodic point $z$  with period $p=\sum_{j=1}^{k} (n_j+m_j+p_j)$ such that
 $$z\in \tilde{B}^{n_1}_{m_1}(x_1,\theta):=\cap_{j=-m_1}^{n_1}f^{-j}B(f^j(x_1),\theta q(f^j(x_1))^{-2})$$
and for $2\leq i\leq k,$ $$f^{\sum_{j=1}^{i-1} (n_j+p_j)+\sum_{j=2}^{i}m_j}(z)\in \tilde{B}^{n_i}_{m_i}(x_i,\theta):=\cap_{j=-m_i}^{n_i}f^{-j}B(f^j(x_i),\theta q(f^j(x_i))^{-2}).$$
\end{Def}

One natural question that arises from Theorem~\ref{Thm:NS} above is

\begin{Que}\label{Que:GenNS}Let $f$ be a $C^{1}$ diffeomorphism. Every $f$ ergodic hyperbolic measure $\mu$ satisfies GNS?
\end{Que}

In particular, if $q(x)\equiv1,$ the required result in Question \ref{Que:GenNS} is obviously valid for uniformly hyperbolic systems \cite{Sigmund}, since  it can be deduced from classical (uniform) specification property \cite{Sigmund} for dynamical balls. Moreover,  $K(\theta,x,m,n)$ can be chosen only dependent on $\theta$ from \cite{Sigmund}.

\bigskip
Here we show that the answer for the Question~\ref{Que:GenNS} is positive, either if $f$ is $C^{1+\alpha}$ or if $f$ is $C^1$ with dominated Oseledec splitting.

\begin{Thm}\label{Thm:GenNS}Let $f$ be a $C^{1+\alpha}\,(\alpha>0)$ diffeomorphism. Then  every $f$ ergodic hyperbolic measure $\mu$ satisfies GNS.

%  That is,  for $\mu$ almost every $x$, any small $\eta>0,$ any $\eta-$slowing varying positive function $q$ (that is, $q(f^{\pm1}(x))\leq e^\eta q(x)$),   any integer $m,\,n,$ and any $\theta>0$ there exists $K:=K(\eta,\theta,x,m,n)$ satisfying
% $$\lim_{\eta\rightarrow 0}\limsup_{m,n\rightarrow+\infty}\frac{K(\eta,\theta,x,m,n)}{m+n}=0$$ and so that the following holds: given points $x_1,x_2,\cdots,x_k$ in a full $\mu$-measure set and positive integers $m_1,\cdots,m_k,n_1,\cdots,n_k,$ there is $ p_i\geq 0 $ with $$\sum_{i=1}^{k} p_i\leq \sum_{i=1}^{k} K(\eta,\theta,x_i,m_i,n_i)$$
% (in particular if $\mu$ is mixing, for any  $t\geq \sum_{i=1}^{k} K(\eta,\theta,x_i,m_i,n_i)$ there is $p_i$ with $\sum_{i=1}^{k} p_i=t$) and a periodic point $z$  with period $p=\sum_{j=1}^{k} (n_j+m_j+p_j)$ such that
%  $$z\in \tilde{B}^{n_1}_{m_1}(x_1,\theta):=\cap_{j=-m_1}^{n_1}f^{-j}B(f^j(x_1),\theta q(f^j(x_1))^{-2})$$
% and for $2\leq i\leq k,$ $$f^{\sum_{j=1}^{i-1} (n_j+p_j)+\sum_{j=2}^{i}m_j}(z)\in \tilde{B}^{n_i}_{m_i}(x_i,\theta):=\cap_{j=-m_i}^{n_i}f^{-j}B(f^j(x_i),\theta q(f^j(x_i))^{-2}).$$
\end{Thm}
\medskip

\begin{Rem}\label{Rem-BowenBall-GNS} In particular, if $q(x)\equiv1,$ we also have a   generalization of Remark \ref{Rem-BowenBall} for
 %dynamical balls as follows. Note that if the function
 %$q(\cdot)\equiv 1$, then the non-uniform dynamical ball  $\widetilde{B}_m^n(x,\theta,q)$ is
 the general well-known dynamical ball $B_m^n(x,\theta)$. I.e., in the conclusion of  GNS we can replace non-uniform dynamical balls by dynamical balls and $\eta$ can be omitted.
\end{Rem}

\bigskip
At the end of this section, we point out that the above results are also valid for  $C^1$ non-uniformly hyperbolic systems with dominated splitting which is based on a recent result of \cite{Tian}.  Before that we recall the notion of dominated splitting. Let $\Delta$ be an $f-$invariant set
and $T_{\Delta}M=E\oplus F$ be a $Df-$invariant splitting on
$\Delta$. $T_{\Delta}M=E\oplus F$ is called $(S_0, \lambda)
$-dominated on $\Delta$ (or simply dominated), if there exist two
constants $S_0\in \mathbb{Z}^+$ and $\lambda>0$ such that
$$\frac1S\log\frac
{\|Df^S|_{E(x)}\|}{m(Df^S|_{F(x)})}\leq-2\lambda,\,\,\forall x \in
\Delta,\,\,S\geq S_0.$$

\begin{Thm}\label{Thm:HypMeaDomSpl} Let $f$ be a $C^{1}$ diffeomorphism. Then every ergodic hyperbolic measure $\mu$ in whose Oseledec splitting the stable bundle dominates the unstable bundle on $supp(\mu)$ satisfies GNS.
\end{Thm}

\section{Pesin theory}\setlength{\parindent}{2em} In this section we give a quick review
concerning some notions and results of $C^{1+\alpha}$ Pesin theory. We point the
reader to \cite{Katok1,Katok2,Pollicott} for more details. Let $f: M\to M$ be a $C^{1+\alpha}$ diffeomorphism. We recall the concept of Pesin set and recall Katok's shadowing lemma in this section.

\smallskip

\subsection{\bf Pesin set}
 Given $\lambda,\mu \gg \varepsilon>0$, and for all $k \in \mathbb{Z}^{+}$,
we define $\Lambda_{k}=\Lambda_{k}(\lambda,\mu;\varepsilon)$ to be
all points $x \in M$ for which there is a splitting
$T_{x}M=E_{x}^{s} \oplus E_{x}^{u}$ with invariant property
$D_{x}f^{m}(E_{x}^{s})=E_{f^{m}x}^{s}$ and
$D_{x}f^{m}(E_{x}^{u})=E_{f^{m}x}^{u}$ satisfying:

\begin{enumerate}
\item[(a)] $\|Df^{n}|_{E_{f^{m}x}^{s}}\| \leq e^{\varepsilon
k}e^{-(\lambda-\varepsilon)n}e^{\varepsilon \mid m\mid},~\forall
m\in\mathbb{Z},~n\geq 1$;

\item[(b)] $\|Df^{-n}|_{E_{f^{m}x}^{u}}\| \leq e^{\varepsilon
k}e^{-(\mu-\varepsilon)n}e^{\varepsilon \mid m\mid},~\forall m\in
\mathbb{Z}, ~n\geq 1$;

\item[(c)] $\tan (\angle(E_{f^{m}x}^{s},E_{f^{m}x}^{u})) \geq
e^{-\varepsilon k}e^{-\varepsilon \mid m\mid},~\forall m\in
\mathbb{Z}$.
\end{enumerate} We set
$\Lambda=\Lambda(\lambda,\mu;\varepsilon)=\bigcup_{k=1}^{+\infty}
\Lambda_{k}$ and call $\Lambda$ a Pesin set.

According to Oseledec Theorem\cite{Os}, every ergodic hyperbolic measure $\mu$ has $s \,\,(s\leq
d=dim M)$ nonzero  Lyapunov exponents
$$\lambda_{1}< \cdot\cdot\cdot<\lambda_{r}<0<\lambda_{r+1}<\cdot\cdot\cdot< \lambda_{s}$$
with associated Oseledec splitting
$$T_{x}M=E_{x}^{1}\oplus\cdot\cdot\cdot\oplus E_{x}^{s},\,\,\,\,x\in O(\mu),$$
where we recall that $O(\mu)$ denotes an Oseledec basin of $\mu.$ If
we denote by $\lambda$ the absolute value of the largest negative
Lyapunov exponent $\lambda_r$ and $\mu$ the smallest positive
Lyapunov exponent $\lambda_{r+1}$, then for any $0<\varepsilon<\min\{\lambda,\,\mu\},$ one has $\mu$ full-measure Pesin set
$\Lambda=\Lambda(\lambda,\mu;\varepsilon)$ (see, for example, Proposition 4.2 in \cite{Pollicott}). And for any point $x\in O(\mu) \cap \Lambda,$ $E^s_x$ and $E^u_x$ coincide with $
E^1_x\oplus\cdot\cdot\cdot\oplus E^r_x$ and $E^{r+1}_x\oplus
\cdot\cdot\cdot\oplus E^s_x$ respectively.

The following statements are elementary properties of Pesin blocks(see \cite{Pollicott}):
\begin{enumerate}
\item[(a)] $\Lambda_{1} \subseteq \Lambda_{2} \subseteq
\Lambda_{3}\subseteq \cdot\cdot\cdot$;

\item[(b)] $f(\Lambda_{k}) \subseteq
\Lambda_{k+1},~f^{-1}(\Lambda_{k}) \subseteq \Lambda_{k+1}$;

\item[(c)] $\Lambda_{k}$ is compact for $\forall\,\, k\geq 1$;

\item[(d)] for $\forall\,\, k\geq 1$ £¬the splitting $x\to E_{x}^{u}\oplus
E_{x}^{s}$ depends continuously on  $x\in\Lambda_{k}$.
\end{enumerate}

\subsection {\bf Shadowing lemma}

We recall Katok's shadowing lemma\cite{Pollicott} in this subsection. Let $(\delta_{k})_{k=1}^{+
\infty}$ be a sequence of positive real numbers. Let
$(x_{n})_{n=-\infty}^{+ \infty}$ be a sequence of points in
$\Lambda=\Lambda(\lambda, \mu, \varepsilon)$ for which there exists
a sequence $(s_{n})_{n=-\infty}^{+ \infty}$ of positive integers
satisfying:

\begin{enumerate}
\item[(a)] $x_{n}\in \Lambda _{s_{n}}, ~\forall n\in \mathbb{Z}$;

\item[(b)] $\mid s_{n}-s_{n-1}\mid \leq 1, ~\forall n\in \mathbb{Z}$;

\item[(c)] $d(fx_{n},x_{n+1})\leq \delta_{s_{n}}, ~\forall n\in \mathbb{Z}$;
\end{enumerate}
then we call $(x_{n})_{n=-\infty}^{+ \infty}$ a
$(\delta_{k})_{k=1}^{+ \infty}$  pseudo-orbit. Given $\theta>0$, a
point $x\in M$ is an
 $\tau$-shadowing point for the $(\delta_{k})_{k=1}^{+ \infty}$
 pseudo-orbit if $d(f^{n}x,x_{n+1})\leq \tau \varepsilon_{s_{n}},~ \forall n\in
\mathbb{Z}$, where $\varepsilon_{k}=\varepsilon_{0}e^{-\varepsilon
k}$ and $\varepsilon_0$ is a constant only dependent on the system of $f$.
\bigskip

\begin{Lem}\label{LemShadow}
(Shadowing lemma) Let $f:M\rightarrow M$ be a
$C^{1+\alpha}$ diffeomorphism, with a non-empty Pesin set
$\Lambda=\Lambda(\lambda,\mu;\varepsilon)$ and fixed parameters,
$\lambda,\mu\gg\varepsilon>0$. For $\forall \tau >0$ there exists a
sequence $(\delta_{k})_{k=1}^{+ \infty}$ such that for any
$(\delta_{k})_{k=1}^{+ \infty}$ pseudo-orbit there exists a unique
$\tau$-shadowing
point.
\end{Lem}

%\begin{Lem}[Closing lemma \cite{Katok2}]\label{closing lem}
%Let $f : M \rightarrow M $ be a $C^{1+\alpha}$ diffeomorphism and let $\Lambda=\Lambda(\alpha,\beta;\,\epsilon)$
%be a nonempty Pesin set. For $\forall\,k\geq1$ and  $\tau> 0$, there exists $b=b(k,\tau)>0$ such that:
%if $x,\,  f^p(x) \in \Lambda_k$ and $d(x,f^px) <b$ then there exists a periodic point $z\in M$ with $z = f^pz$ and
%$$d(f^i(z),f^i(x))<\tau \varepsilon_0 e^{\min\{i,p-i\}\varepsilon},\,i=0,1,2,\cdots,p.$$\end{Lem}

\section{Recurrence Times}

In this section we always assume that $f:X\rightarrow X$ is a
homeomorphism on a compact metric space, $\mu$ is an invariant measure  and $\Gamma$ is an subset of $X$ with positive measure for $\mu$. For $x\in \Gamma,$ define
 $$
\cdots<t_{-2}(x)<t_{-1}(x)<t_0(x)=0<t_1(x)<t_2(x)<\cdots
$$
to be the all times
such that $f^{t_i}(x)\in\Gamma$ (called recurrence times). By Poincar\'{e} Recurrence Theorem, the sequence above is well-defined at $\mu$ a.e. $x\in \Gamma.$ Note that
 $$
t_1(f^{t_i}(x))=t_{i+1}(x)-t_i(x).
$$
In general, we call $t_1$ to be the first recurrence time. An increasing sequence of natural numbers $\{N_i\}_{i\geq 1}$ is called {\it nonlacunary}, if $\lim_{i\rightarrow+\infty}\frac{N_{i+1}}{N_i}=1 $.
For recurrence times, we have a basic   proposition of their nonlacunary as follows.

\begin{Prop}\label{Prop:Rec2} For $\mu$ a.e $x\in \Gamma,$
$$\lim_{i\rightarrow+\infty}\frac{t_{i+1}(x)}{t_i(x)}=1\,\,\,\,\text{and}\,\,\, \lim_{i\rightarrow-\infty}\frac{t_{i-1}(x)}{t_i(x)}=1.$$
In particular, $$\lim_{i,j\rightarrow+\infty}\frac{t_{i+1}(x)-t_{-j-1}(x)}{t_i(x)-t_{-j}(x)}=1.$$
\end{Prop}

{\bf Proof} It can be proved by using  Borel-Cantelli lemma and Ka${\check{c}}$ Lemma( See \cite{OV08}, Proposition 3.8). Here we give another direct proof of the first equality which only depends on Birkhoff Ergodic Theorem.
The proof of the remain equalities are similar.

By Birkhoff Ergodic Theorem, for any subset $A\subseteq M,$ the limit function $$\chi^*_A(x):=\lim_{n\rightarrow+\infty}\frac1n\sum_{j=0}^{n-1}\chi_A(f^j(x))$$ exists for $\mu$ a.e. $x$ and $f-$invariant. Moreover, $\int\chi^*_A(x)d\mu=\int\chi_A(x)d\mu=\mu(A).$

If $A:=\Gamma,$ we claim that for $\mu$ a.e. $x\in \Gamma$, $\chi^*_\Gamma(x)>0$ (this is a basic fact of recurrence and we will give the proof below). Then, by the definition of $t_i(x)$, for $\mu$ a.e $x\in \Gamma,$ $$\lim_{i\rightarrow+\infty}\frac{i}{t_i(x)}=\lim_{t_i(x)\rightarrow+\infty}\frac1{t_i(x)}\sum_{j=0}^{t_i(x)-1}\chi_\Gamma(f^j(x))=\chi^*_\Gamma(x)>0.$$
Thus, for $\mu$ a.e $x\in \Gamma,$
$$\lim_{i\rightarrow+\infty}\frac{t_{i+1}(x)}{t_i(x)}=\lim_{i\rightarrow+\infty}\frac{t_{i+1}(x)}{i+1}\cdot\frac{i+1}{i}\cdot\frac{i}{t_i(x)}
=\frac1{\chi^*_\Gamma(x)}\cdot 1 \cdot\chi^*_\Gamma(x)=1.$$

Now we start to prove the claim. If $\mu$ is ergodic, it is obvious  since $\chi^*_\Gamma(x)\equiv\int\chi_\Gamma(x)d\mu=\mu(\Gamma)>0$ holds for a.e. point $x$. For the general invariant case, we prove by contradiction. Assume that the set $\Gamma_1:=\{x\in\Gamma\,|\,\chi^*_\Gamma(x)=0 \}$ has $\mu$ positive measure. Consider $A:=\Gamma_1,$ then $\chi^*_{\Gamma_1}(x)$ exists for $\mu$ a.e. $x$ and $\int\chi^*_{\Gamma_1}(x)d\mu=\int\chi_{\Gamma_1}(x)d\mu=\mu(\Gamma_1)>0.$ Thus there exists $x \in \bigcup_{n\in\mathbb{Z}}f^n(\Gamma_1)$ such that $\chi^*_{\Gamma_1}(x)>0,$ since by definition  $\chi^*_{\Gamma_1}(x)\equiv0$ for all $x\in (\bigcup_{n\in\mathbb{Z}}f^n(\Gamma_1))^c.$ Recall that $\chi^*_{\Gamma_1}(x)$ is $f$-invariant and thus if we let $y\in\Gamma_1$ such that $f^m(y)=x$ for some integer $m,$ then $\chi^*_{\Gamma_1}(y)=\chi^*_{\Gamma_1}(f^{-m}(x))=\chi^*_{\Gamma_1}(x)>0.$ Note that $\Gamma_1\subset \Gamma$ implies $\chi_{\Gamma_1}(\cdot)\leq\chi_{\Gamma}(\cdot)$ and thus $\chi^*_{\Gamma_1}(\cdot)\leq\chi^*_{\Gamma}(\cdot).$ So $$\chi^*_{\Gamma}(y)\geq\chi^*_{\Gamma_1}(y)>0.$$ This contradicts the choice of $\Gamma_1$ since $y\in\Gamma_1$.\qed
\bigskip
\begin{Rem}\label{Rem-RecTime1}  We emphasize the main used technique in Proposition \ref{Prop:Rec2} that for an increasing sequence of integers $\{N_i\}_{i\geq 1}$, one sufficient condition to get nonlacunary is that $\lim_{i\rightarrow+\infty}\frac{i}{N_i}>0$. Thus this can also be used to characterize the nonlacunary of hyperbolic times\cite{Oli} instead of using Borel-Cantelli lemma.
\end{Rem}
\medskip

\begin{Rem}\label{Rem-RecTime} We point out another equivalent statement of Proposition \ref{Prop:Rec2}. That is, $x$ satisfies $\lim_{i\rightarrow+\infty}\frac{t_{i+1}(x)}{t_i(x)}=1$ $\Leftrightarrow$ for any $\epsilon>0,$ there exists large integer $N(x)$ such that for all $n\geq N(x),$ there is $t\in[n,n+n\epsilon)$ such that $t=t_i(x)$ for some $i.$  The proof of this relation is trivial but the technique of the later version is also useful and has been essentially used in the proof to  establish an (in)equality between metric entropy of hyperbolic ergodic measure and the number of hyperbolic periodic points in \cite{Katok2,BP,LST1}. One main technique used in in \cite{Katok2,BP,LST1} is that, the recurrence time varies so slowly (linearly) that if the cardinality of a sequence of sets(for example, separated sets) is growing exponentially, then there exists a new sequence composed of their subsets such that the cardinality still grows exponentially and in every subset the recurrence times   are same for all points.
\end{Rem}

\section{Proof of Theorem \ref{Thm:NS}: NS}
 Set $\tilde \Lambda_k=supp(
\mu|_{\Lambda_{k}})$ and $\tilde
\Lambda=\cup_{k=1}^\infty\tilde \Lambda_k.$ Clearly, $f^{\pm
1} \big(\tilde \Lambda_k \big) \subset \tilde \Lambda_{k +1},$ and the
sub-bundles $E^s(x),\,\, E^u(x)$ depend continuously on $x\in\tilde
\Lambda_{k}.$  Moreover,  $\tilde \Lambda$ is $f-$invariant with
$\mu $-full measure.  Let $\Delta_k\subseteq\tilde \Lambda_k$ be the set of all points of $x$ satisfying that

\textit{(i)}  recurrence times of $x$ are well defined for the set $\Gamma:=\tilde \Lambda_k$ and
\begin {equation}\label{Eq:0Rec}(ii)\,\,\,\,\,\,\,\,\,\,\,\,\,\,\,\,\,\,\,\,\,\,\,\,\,\,\,
\,\,\,\,\,\,\,\,\,\,\,\,\,\,
\lim_{i,j\rightarrow+\infty}\frac{t_{i+1}(x)-t_{-j-1}(x)}{t_i(x)-t_{-j}(x)}=1.\,\,\,\,\,\,\,\,\,\,\,\,\,\,\,\,\,\,\,\,\,\,\,\,\,\,\,\,\,\,\,\,\,\,\,\,\,\,\,\,\,\,\,\,\,\,\,\,\,\,\,\,\,\,\,\,\,\,\,\,\,\,\,\,\,\,
\end {equation}
By Poincar\'{e} Recurrence Theorem and Proposition \ref{Prop:Rec2}, $\mu(\Delta_k)=\mu(\tilde \Lambda_k)$ and thus $\cup_{k\geq1}\Delta_k$ is also a subset of $M$ with full measure. So we only need to prove that for every fixed $\Delta_k$ with positive measure, all points in $\Delta_k$ satisfy the conditions of non-uniform specification property.
Take and fix a point $x\in\Delta_k$.

Recall $\varepsilon$ to be the number that appeared in the definition of Pesin set. Let $0<\eta\leq \varepsilon/2$ and $q$ be an $\eta-$slowing varying positive function, $\theta>0$ and let $m,\,n$ be two positive integers.
Let $\tau=\frac{\theta q^{-2}(x)}{\varepsilon_0}
>0$. By Lemma \ref{LemShadow} for this $\tau$ there exists a sequence
$(\delta_{k})_{k=1}^{+ \infty}$ such that for any
$(\delta_{k})_{k=1}^{+ \infty}$ pseudo-orbit there exists a unique
$\tau$-shadowing point.

Take and fix for $\tilde\Lambda_k$
a finite cover $\alpha_k=\{V_1,V_2,\cdots,V_{r_k}\}$ by
nonempty open balls $V_i$ in $M$ such that diam$(U_i)
<\delta_{k+1}$ and $\mu(U_i)>0$ where $U_i=V_i\cap
\tilde\Lambda_k,\,i=1,\,2,\,\cdots,\,r_k$. This is obtained from the definition of $\tilde\Lambda_k$.  By Birkhoff ergodic theorem and the ergodicity of  $\mu$ we have
\begin {equation}\label{Eq:1}
\lim_{n\rightarrow +\infty}\frac 1n
  \sum_{h=0}^{n-1}\mu(f^{-h}(U_i)\cap U_j)=\mu(U_i)\mu(U_j)>0.
\end{equation}
 Then one can
take
$$ X_{i,\,j}:=\min\{h\in \mathbb{N}\,\,|\,\, h\geq
 1, \,\,\mu(f^{-h}(U_i)\cap
  U_j)>0\}.
$$
By  (\ref{Eq:1}), $1\leq X_{i,\,j}<+\infty$. Let
$$M_{k}=
\max_{1\leq i, j\leq r_k}X_{i,\,j}.$$
Note that $M_{k}$ is a positive integer dependent on $k,\theta$ and the $\eta-$slowing varying positive function $q,$ but independent of $m,n.$ So \begin {equation}\label{Eq:2M_k} \frac{M_k}m,\,\,\,\frac{M_k}n\rightarrow 0
\end {equation} as $m,n\rightarrow \infty$ respectively.

Now we start to find the needed shadowing periodic orbit. Take positive integers $l_1$ and $l_2$  such that
\begin {equation}\label{Eq:2}
 t_{-l_1}< -m \leq  t_{-l_1+1}\,\,\,\,\text{and}\,\,\,\, t_{l_2}>n \geq  t_{l_2-1},
\end{equation}
and take positive integers $s_1$ and $s_2$  such that
 \begin {equation}\label{Eq:3}
  t_{-l_1-s_1}\leq (1+ \frac{2\eta} \varepsilon)  t_{-l_1}< t_{-l_1-s_1+1}\,\,\,\,\text{and}\,\,\,\, t_{l_2+s_2}\geq (1+ \frac{2\eta} \varepsilon)  t_{l_2}> t_{l_2+s_2-1}.
\end{equation}
Since $f^{ t_{-l_1-s_1}}(x), f^{ t_{l_2+s_2}}(x)\in\tilde\Lambda_k,$ we can take $U_i$ and $U_j$ such that
$$f^{ t_{-l_1-s_1}}(x)\in U_i, f^{ t_{l_2+s_2}}(x)\in U_j.$$ By  (\ref{Eq:1}), there exist $y\in U_j$ and $0\leq N\leq M_k$ such that $f^N(y)\in U_i.$

Recall the property of Pesin blocks that $f^{\pm}(\Lambda_k)\subseteq \Lambda_{k+1}.$ Thus, if $u\in\Lambda_k$, then $f^{i}(u)\in \Lambda_{k+|i|},\,\forall\,i\in\mathbb{Z}.$
Note that $$f^{ t_{-l_1-s_1}}(x),x, f^{ t_{l_2+s_2}}(x), y, f^N(y)\in \tilde\Lambda_k\subseteq \Lambda_k$$ and $$d(f^{ t_{l_2+s_2}}(x), y)<\delta_{k+1},\,\,d(f^{ t_{-l_1-s_1}}(x), f^N(y))<\delta_{k+1}.$$ So $$ f^{ t_{-l_1-s_1}}(x)\in  \Lambda_k, f^{ t_{-l_1-s_1+1}}(x)\in  \Lambda_{k+1},\cdots f^{ t_{-l_1-s_1+i}}(x)\in  \Lambda_{\min\{k+i,k+l_1+s_1-i\}}, $$ $$\cdots,f^{-1}(x)\in  \Lambda_{k+1},x\in  \Lambda_{k}, f(x)\in  \Lambda_{k+1},\cdots, f^i(x)\in  \Lambda_{\min\{k+i,k+l_2+s_2-i\}},$$$$\cdots, f^{ t_{l_2+s_2-1}}(x)\in  \Lambda_{k+1}, y\in  \Lambda_k,\cdots,f^{i-1}(y)\in  \Lambda_{\min\{k+i,k+N-i\}},\cdots, f^{N-1}(y)\in  \Lambda_{k+1}.$$ Repeat the above sequence of points infinitely many times and thus we get a
$(\delta_{k})_{k=1}^{+ \infty}$ pseudo-orbit. Then there exists a unique
$\tau$-shadowing point $z$. Note that $f^{ t_{l_2+s_2}- t_{-l_1-s_1}+N}(z)$ is also a $\tau$-shadowing point. So $$f^{p}(z)=z$$ where $p= t_{l_2+s_2}- t_{-l_1-s_1}+N.$

 Now we start to verify the conditions of non-uniform specification property for the above chosen shadowing point $z$. Define $K:=K(\eta,\theta,x,m,n)= t_{l_2+s_2}- t_{-l_1-s_1}+M_k-m-n.$ Clearly we have $p\leq m+n+K$ since $N\leq M_k,$  and thus for the first condition of NS we only need to show $$z\in\tilde{B}^n_m(x,\theta):=\cap_{i=-m}^nf^{-i}B(f^k(x),\theta q(f^k(x))^{-2}).$$  Firstly, we consider $0\leq i\leq  t_{l_2}-1$ and calculate $d(f^i(x),f^i(z)).$ More precisely, $\tau-$shadowing implies that
\begin{eqnarray*}
& &d(f^i(x),f^i(z))\nonumber\\
&\leq& \max\{\tau \varepsilon_{k+i}, \tau \varepsilon_{k+ t_{l_2+s_2}-i}\}\nonumber\\
&\leq& \max\{\tau \varepsilon_i, \tau \varepsilon_{ t_{l_2+s_2}- t_{l_2}}\}\,\,\,\,\,\,\,(\text{note that } i\leq t_{l_2} \text{ and } \varepsilon_k \text{ is a decreasing sequence})\nonumber\\
&=& \max\{\tau \varepsilon_0e^{-i\varepsilon}, \tau \varepsilon_0e^{-( t_{l_2+s_2}- t_{l_2})\varepsilon}\}\nonumber\\
&\leq& \max\{\tau \varepsilon_0e^{-i\varepsilon}, \tau \varepsilon_0e^{-2 t_{l_2}\eta}\}\,\,\,\,\,\,\,(\text{using } (\ref{Eq:3}))\nonumber\\
&\leq& \max\{\tau \varepsilon_0e^{-2i\eta}, \tau \varepsilon_0e^{-2i\eta}\}\,\,\,\,\,\,\,(\text{using } 2\eta<\varepsilon\text{ and }i\leq  t_{l_2} )\nonumber\\
&=& \tau \varepsilon_0e^{-2i\eta}\nonumber\\
&=& \theta q^{-2}(x)e^{-2i\eta}\,\,\,\,\,\,\,(\text{by the choice of } \tau )\nonumber\\
&\leq& \theta q(f^i(x))^{-2}\,\,\,\,\,\,\,\,(\text{using }q(f^i(x))\leq q(x)e^{i\eta}).\nonumber
\end{eqnarray*}
Secondly we can follow the similar method to show that for $ t_{-l_1}+1\leq i\leq 0,$
\begin{eqnarray*}
d(f^i(x),f^i(z))\leq \theta q(f^i(x))^{-2}.\nonumber
\end{eqnarray*}
Notice that $ t_{-l_1}<-m$ and $n<  t_{l_2}$ and thus $z\in\tilde{B}^n_m(x,\theta)$.

At the end  we prove   the second condition of NS: $$\lim_{\eta\rightarrow 0}\limsup_{m,n\rightarrow+\infty}\frac{K(\eta,\theta,x,m,n)}{m+n}=0.$$ In fact, by using (\ref{Eq:0Rec}), (\ref{Eq:2M_k}), (\ref{Eq:2}) and (\ref{Eq:3}), we have
\begin{eqnarray*}
& &\limsup_{m,n\rightarrow+\infty}\frac{K(\eta,\theta,x,m,n)}{m+n}\nonumber\\
&\leq&\limsup_{n\rightarrow+\infty}\frac{ t_{l_2+s_2}- t_{-l_1-s_1}}{m+n}
+\limsup_{m,n\rightarrow+\infty}\frac{M_k}{m+n}-1\nonumber\\
&\leq&\limsup_{m,n\rightarrow+\infty}\frac{ t_{l_2+s_2}- t_{-l_1-s_1}}{ t_{l_2-1}- t_{-l_1+1}}
+0-1\,\,\,\,\,\,\,\,(\text{using}\,\,  (\ref{Eq:2}), (\ref{Eq:2M_k}) )\nonumber\\
&=&\limsup_{m,n\rightarrow+\infty}(\frac{ t_{l_2+s_2}- t_{-l_1-s_1}}{ t_{l_2+s_2-1}- t_{-l_1-s_1+1}}
\cdot\frac{ t_{l_2+s_2-1}- t_{-l_1-s_1+1}}{ t_{l_2}- t_{-l_1}}\cdot\frac{ t_{l_2}- t_{-l_1}}{ t_{l_2-1}- t_{-l_1+1}})
-1\,\,\,\nonumber\\
&=&\limsup_{m,n\rightarrow+\infty}\frac{ t_{l_2+s_2-1}- t_{-l_1-s_1+1}}{ t_{l_2}- t_{-l_1}}
-1\,\,\,\,\,\,\,\,(\text{using}\,\,  (\ref{Eq:0Rec}) )\nonumber\\
&\leq&1+\frac{2\eta}\varepsilon
-1=\frac{2\eta}\varepsilon\,\,\,\,\,\,\,(\text{using}\,\,  (\ref{Eq:3}) ).\nonumber
\end{eqnarray*}
Letting $\eta\rightarrow 0$, one has
$$
\lim_{\eta\rightarrow 0}\limsup_{m,n\rightarrow+\infty}\frac{K(\eta,\theta,x,m,n)}{m+n}=0.
$$
So we complete the proof.\qed

\bigskip

\begin{Rem}\label{Rem-NS}  In particular, if  $\mu$ is a mixing hyperbolic measure, we can replace inequality (\ref{Eq:1}) by \begin {equation}\label{Eq:Remark}
\lim_{n\rightarrow +\infty}\mu(f^{-n}(U_i)\cap U_j)=\mu(U_i)\mu(U_j)>0.
\end{equation}
 Then  by  (\ref{Eq:Remark}) we can
take a finite integer
$$ X_{i,\,j}=\max\{n\in \mathbb{N}\,\,|\,\, n\geq
 1, \,\,\mu(f^{-n}(U_i)\cap
  U_j)=0\}+1.
$$ Let
$$M_{k}=
\max_{1\leq i, j\leq r_k}X_{i,\,j}.$$
Then for any $ N\geq M_k$ there exist $y\in U_j$  such that $f^N(y)\in U_i.$
So we can follow the above proof and then the non-uniform specification can be stronger: for any $p\geq n+m+K,$ the non-uniform dynamical ball $\tilde{B}^n_m(x,\theta)$ contains a periodic point with period $p$.

\end{Rem}

\section{Proof of Theorem~\ref{Thm:GenNS}: GNS%Generalized Non-uniform Specification
}

In this section we prove Theorem \ref{Thm:GenNS}. Before, we show two  propositions as follows.

\begin{Prop}\label{Prop:ExpGenNS}Let $f$ be a $C^{1+\alpha}\,(\alpha>0)$ diffeomorphism. Then for any small $0<\sigma<1$, there is a subset $\Lambda^*_\sigma$ with $\mu(\Lambda^*_\sigma)>1-\sigma$ such that
for every $x\in \Lambda^*_\sigma$,  any small $\eta>0,$  any $\theta_*>0$  and  any integer $m,\,n,$ there exists $K:=K(\eta,\theta_*,x,m,n)$ satisfying
$$\lim_{\eta\rightarrow 0}\limsup_{m,n\rightarrow+\infty}\frac{K(\eta,\theta_*,x,m,n)}{m+n}=0$$ and so that the following holds: given points $x_1,x_2,\cdots,x_k$ in $\Lambda^*_\sigma$ and positive integers $m_1,\cdots,m_k,n_1,\cdots,n_k,$ there are numbers $ p_{*i}\geq 0\,(1\leq i\leq k)$ with $$\sum_{i=1}^{k} p_{*i}\leq \sum_{i=1}^{k} K(\eta,\theta_*,x_i,m_i,n_i)$$(in particular if $\mu$ is mixing, for any  $t\geq \sum_{i=1}^{k} K(\eta,\theta_*,x_i,m_i,n_i)$ there is $p_{*i}$ with $\sum_{i=1}^{k} p_{*i}=t$) and there exists a periodic point $z_*$  with period $p_*=\sum_{j=1}^{k} (n_j+m_j+p_{*j})$ such that
 $$z_*\in \cap_{j=-m_1}^{n_1}f^{-j}\big(B(f^j(x_1),\theta_* e^{-2|j|\eta})\big)$$
and for $2\leq i\leq k,$ $$f^{\sum_{j=1}^{i-1} (n_j+p_{*j})+\sum_{j=2}^{i}m_j}(z_*)\in \cap_{j=-m_i}^{n_i}f^{-j}B(f^j(x_i),\theta_* e^{-2|j|\eta}).$$

This property is also valid for one side case.
\end{Prop}

{\bf Proof}
The proof is a generalization of that of Theorem \ref{Thm:NS}. Recall $\varepsilon$ to be the number that appeared in the definition of Pesin set. Recall $\Delta_{k_*}$ to be the set introduced in the proof of Theorem \ref{Thm:NS}
and note that if ${k_*}$ is large then the measure of $\Delta_{k_*}$ is close to 1. Let $k_*$ be large enough such that  $\Delta_{k_*}$
satisfies $\mu(\Delta_{k_*})>1-\sigma$. We will prove this $\Delta_{k_*}$ is the required $\Lambda^*_\sigma$.

 Let $\tau=\frac{\theta_*}{\varepsilon_0}>0$.
By Lemma \ref{LemShadow} there exists a sequence
$(\delta_{k})_{k=1}^{+ \infty}$ such that for any
$()_{k=1}^{+ \infty}$ pseudo-orbit there exists a unique
$\tau$-shadowing point.

Let $x\in\Delta_{k_*}$,  $0<\eta\leq \varepsilon/2$ and $m,\,n$ be two positive integers. Note that $\delta_{k_*}$ dependent on $\tau=\frac{\theta_*}{\varepsilon_0}$ and $\Lambda_{k_*}$, but independent of $x$. This is different to the one in the proof of Theorem \ref{Thm:NS}. Let $K(\eta,\theta_*,x,m,n)$ be the number defined as in the proof of Theorem \ref{Thm:NS}.(Note that the choices of $l_1,s_1$ and $l_2,s_2$ only depends on $\Delta_{k_*}$, $x,m,n,$ and the choice of $M_{k_*}$ only depends on $\delta_{{k_*}+1}$ and thus $K(\eta,\theta_*,x,m,n)$ only depends $\Delta_{k_*}$ and $\tau=\frac{\theta_*}{\varepsilon_0}$). So  this $K(\eta,\theta_*,x,m,n)$ also satisfies $$\lim_{\eta\rightarrow 0}\limsup_{m,n\rightarrow+\infty}\frac{K(\eta,\theta_*,x,m,n)}{m+n}=0.$$
Re-denote the $l_1,s_1$ and $l_2,s_2$ with respect to $x$ in the proof of Theorem \ref{Thm:NS} by $l_1(x),s_1(x)$ and $l_2(x),s_2(x)$.

Given points $x_1,x_2,\cdots,x_k$ in $\Delta_{k_*}$ and positive integers $m_1,\cdots,m_k,n_1,\cdots,n_k,$
similar as the proof of Theorem \ref{Thm:NS}, we can take $$ y_i, f^{N_i}(y_i)\in \tilde\Lambda_{k_*}$$ with $0\leq N_i\leq M_{k_*}$ such that $$d(f^{t_{l_2(x_i)+s_2(x_i)}}(x_i), y_i)<\delta_{{k_*}+1},\,\,d(f^{t_{-l_1(x_{i+1})-s_1(x_{i+1})}}(x_{i+1}), f^{N_i}(y_i))<\delta_{{k_*}+1},\,\,\,\forall \,\,1\leq i\leq {k_*},$$ where $x_{k+1}=x_1.$ Note that $$f^{t_{-l_1-s_1}}(x_i),x, f^{t_{l_2+s_2}}(x_i), y_i, f^{N_i}(y_i)\in \tilde\Lambda_{k_*}\subseteq \Lambda_{k_*}.$$ Similar as the proof of Theorem \ref{Thm:NS}, by shadowing lemma there is a periodic point $z_*$  with period $p_*=\sum_{j=1}^{k} (t_{l_2(x_j)+s_2(x_j)}-t_{l_1(x_j)+s_1(x_j)}+N_j)$ such that
 $$z_*\in \cap_{j=-m_1}^{n_1}f^{-j}B(f^j(x_1),\tau\varepsilon_0  e^{-2|j|\eta})=\cap_{j=-m_1}^{n_1}f^{-j}B(f^j(x_1),\theta_*  e^{-2|j|\eta})$$
and for $2\leq i\leq k,$ $$f^{\sum_{j=1}^{i-1} (t_{l_2(x_j)+s_2(x_j)}+N_j)-\sum_{j=2}^{i}t_{l_1(x_j)+s_1(x_j)}}(z_*)\in \cap_{j=-m_i}^{n_i}f^{-j}B(f^j(x_i),\theta_*  e^{-2|j|\eta}).$$

Let $ p_{*i}=t_{l_2(x_i)+s_2(x_i)}-n_i+N_i-t_{-l_1(x_{i+1})-s_1(x_{i+1})}-m_{i+1} $ where $m_{k+1}=m_1$ and thus $\sum_{i=1}^{k} p_{*i}=\sum_{i=1}^{k}(t_{l_2(x_i)+s_2(x_i)}-t_{-l_1(x_{i})-s_1(x_i)}+N_i-m_{i}-n_i)\leq \sum_{i=1}^{k} K(\eta,\theta_*,x_i,m_i,n_i)$. Then
the  periodic point $z_*$  satisfies that its period is $p_*=\sum_{j=1}^{k} (n_j+m_j+p_{*j})$,
 $$z_*\in \cap_{j=-m_1}^{n_1}f^{-j}B(f^j(x_1),\theta_*  e^{-2|j|\eta})$$
and for $2\leq i\leq k,$ $$f^{\sum_{j=1}^{i-1} (n_j+p_{*j})+\sum_{j=2}^{i}m_j}(z_*)=f^{\sum_{j=1}^{i-1} (t_{l_2(x_j)+s_2(x_j)}+N_j)-\sum_{j=2}^{i}t_{l_1(x_j)+s_1(x_j)}}(z_*)$$
$$\in \cap_{j=-m_i}^{n_i}f^{-j}B(f^j(x_i),\theta_*  e^{-2|j|\eta}).\,\,\,\,\,\,\,\qed$$

\begin{Prop}\label{Prop:GenNS}Let $f$ be a $C^{1+\alpha}\,(\alpha>0)$ diffeomorphism. Then for any small $0<\sigma<1$, any small $\eta>0,$ any $\eta-$slowing varying positive function $q$ (that is, $q(f^{\pm1}(x))\leq e^\eta q(x)$) and any $\theta>0$, there is a subset $\Lambda_\sigma$ with $\mu(\Lambda_\sigma)>1-\sigma$ such that
for every $x\in \Lambda_\sigma$,    any integer $m,\,n,$ there exists $K_*:=K_*(\eta,\theta,x,m,n)$ satisfying
$$\lim_{\eta\rightarrow 0}\limsup_{m,n\rightarrow+\infty}\frac{K_*(\eta,\theta,x,m,n)}{m+n}=0$$ and so that the following holds: given points $x_1,x_2,\cdots,x_k$ in $\Lambda_\sigma$ and positive integers $m_1,\cdots,m_k,n_1,\cdots,n_k,$ there is $ p_{*i}\geq 0$ with $$\sum_{i=1}^{k} p_{*i}\leq \sum_{i=1}^{k} K_*(\eta,\theta,x_i,m_i,n_i)$$(in particular if $\mu$ is mixing, for any  $t\geq \sum_{i=1}^{k} K_*(\eta,\theta,x_i,m_i,n_i)$ there is $p_{*i}$ with $\sum_{i=1}^{k} p_{*i}=t$) and a periodic point $z_*$  with period $p_*=\sum_{j=1}^{k} (n_j+m_j+p_{*j})$ such that
 $$z_*\in \tilde{B}^{n_1}_{m_1}(x_1,\theta):=\cap_{j=-m_1}^{n_1}f^{-j}B(f^j(x_1),\theta q(f^j(x_1))^{-2})$$
and for $2\leq i\leq k,$ $$f^{\sum_{j=1}^{i-1} (n_j+p_{*j})+\sum_{j=2}^{i}m_j}(z_*)\in \tilde{B}^{n_i}_{m_i}(x_i,\theta):=\cap_{j=-m_i}^{n_i}f^{-j}B(f^j(x_i),\theta q(f^j(x_i))^{-2}).$$

This property is also valid for one side case.
\end{Prop}

{\bf Proof}
Recall $\varepsilon$ to be the number that appeared in the definition of Pesin set. Let $0<\eta\leq \varepsilon/2$ and $q$ be an $\eta-$slowing varying positive function. Recall $\Lambda^*_{\frac{\sigma}2}$ to be the set introduced in the proof of Proposition \ref{Prop:ExpGenNS} for $\frac{\sigma}2$.
Let  $\theta_*>0$ be small enough such that we can take $\Lambda_\sigma\subseteq\Lambda^*_{\frac{\sigma}2}$ with $\mu(\Lambda_\sigma)>1-\sigma$ and every point $x\in\Lambda_\sigma$ satisfies $\theta q^{-2}(x)\geq \theta_*.$

Let $x\in\Lambda_\sigma\subseteq\Lambda^*_{\frac{\sigma}2}$, $m,\,n$ be two positive integers. Let $K(\eta,\theta_*,x,m,n)$ be the number as in Proposition \ref{Prop:ExpGenNS}, then if we can take  $K_*(\eta,\theta,x,m,n):=K(\eta,\theta_*,x,m,n)$ and thus  this  $K_*(\eta,\theta,x,m,n)$ also satisfies $$\lim_{\eta\rightarrow 0}\limsup_{m,n\rightarrow+\infty}\frac{K_*(\eta,\theta,x,m,n)}{m+n}=0.$$
Re-denote the $l_1,s_1$ and $l_2,s_2$ with respect to $x$ in the proof of Theorem \ref{Thm:NS} by $l_1(x),s_1(x)$ and $l_2(x),s_2(x)$.

Given points $x_1,x_2,\cdots,x_k$ in $\Lambda_\sigma\subseteq\Lambda^*_{\frac{\sigma}2}$ and positive integers $m_1,\cdots,m_k,n_1,\cdots,n_k,$
 by Proposition \ref{Prop:ExpGenNS} there is $ p_{*i}\geq 0$ with $$\sum_{i=1}^{k} p_{*i}\leq \sum_{i=1}^{k} K(\eta,\theta_*,x_i,m_i,n_i)$$(in particular if $\mu$ is mixing, for any  $t\geq \sum_{i=1}^{k} K(\eta,\theta_*,x_i,m_i,n_i)$ there is $p_{*i}$ with $\sum_{i=1}^{k} p_{*i}=t$) and there is a periodic point $z_*$  with period $p_*=\sum_{j=1}^{k} (n_j+m_j+p_{*j})$,
 $$z_*\in \cap_{j=-m_1}^{n_1}f^{-j}B(f^j(x_1),\theta_*  e^{-2|j|\eta})$$
and for $2\leq i\leq k,$ $$f^{\sum_{j=1}^{i-1} (n_j+p_{*j})+\sum_{j=2}^{i}m_j}(z_*)\in \cap_{j=-m_i}^{n_i}f^{-j}B(f^j(x_i),\theta_*  e^{-2|j|\eta}).\,\,\,\,\,\,\,$$
Since $x_i\in\Lambda_\sigma$, then $\theta q^{-2}(x_i)\geq \theta_*.$ Using $q(f^{j}(x_i))\leq e^{|j|\eta} q(x_i)$, we have  $$z_*\in \cap_{j=-m_1}^{n_1}f^{-j}B(f^j(x_1),\theta_*  e^{-2|j|\eta})\subseteq \cap_{j=-m_1}^{n_1}f^{-j}B(f^j(x_1),\theta q^{-2}(x_1)  e^{-2|j|\eta})$$
$$\subseteq \cap_{j=-m_1}^{n_1}f^{-j}B(f^j(x_1),\theta q^{-2}(f^j(x_1)) )$$
and similarly for $2\leq i\leq k,$ $$f^{\sum_{j=1}^{i-1} (n_j+p_{*j})+\sum_{j=2}^{i}m_j}(z_*)\in \cap_{j=-m_i}^{n_i}f^{-j}B(f^j(x_i),\theta q^{-2}(f^j(x_i))).\,\,\,\,\,\,\,$$

\bigskip

Now we start to prove Theorem \ref{Thm:GenNS}.

{\bf Proof of  Theorem \ref{Thm:GenNS}} We use Proposition \ref{Prop:GenNS} to give a proof.
 Let   $\Lambda_\sigma$ be a \textit{fixed} positive $\mu$-measure set from Proposition \ref{Prop:GenNS}.
%and we may assume every point $x\in\Lambda_\sigma$ can be recurrent from Poincar\'{e} Recurrence Theorem.
Note that  $\cup_{j\geq 0}f^j(\Lambda_\sigma)$ is a full $\mu$-measure set by the ergodicity of $\mu$. We consider $x\in \cup_{j\geq 0}f^j(\Lambda_\sigma)$. Clearly we can take a finite (and fixed) number $t(x)\geq 0$(which can be chosen the first time) such that $x\in f^{t(x)}(\Lambda_\sigma)$ and thus $f^{-t(x)}(x)\in \Lambda_\sigma.$
Let $K_*(\eta,\theta,f^{-t(x)}(x),m,n+t(x))$ be the number as in Proposition \ref{Prop:GenNS} and define $$K(\eta,\theta,x,m,n):=K_*(\eta,\theta,f^{-t(x)}(x),m,n+t(x))+t(x).$$ Then $K(\eta,\theta,x,m,n)$ satisfies
$$\lim_{\eta\rightarrow 0}\limsup_{m,n\rightarrow+\infty}\frac{K(\eta,\theta,x,m,n)}{m+n}$$$$=\lim_{\eta\rightarrow 0}\limsup_{m,n+t(x)\rightarrow+\infty}\frac{K_*(\eta,\theta,f^{-t(x)}(x),m,n+t(x))+t(x)}{m+n+t(x)}=0.$$

Given $x_1,\cdots,x_k \in \cup_{j\geq 0}f^j(\Lambda_\sigma)$ and positive integers $m_1,m_2,\cdots,m_k,n_1,\cdots,n_k$ large enough,  we consider  points $f^{-t(x_1)}(x_1),\cdots,f^{-t(x_k)}(x_k) \in \Lambda_\sigma$ and positive integers  $$m_1,m_2,\cdots,m_k,n_1+t(x_1), \cdots, n_k+t(x_k),$$ by Proposition \ref{Prop:GenNS} there is
$ p_{*i}\geq 0$ with $$\sum_{i=1}^{k} p_{*i}\leq \sum_{i=1}^{k} K_*(\eta,\theta,f^{-t(x_i)}(x_i),m_i,n_i+t(x_i))$$(in particular if $\mu$ is mixing, for any  $t\geq \sum_{i=1}^{k} K_*(\eta,\theta,f^{-t(x_i)}(x_i),m_i,n_i+t(x_i))$ there is $p_{*i}$ with $\sum_{i=1}^{k} p_{*i}=t$)
and a periodic point $z_*$  with period $$p_*=\sum_{j=1}^{k} (m_j+n_j+t(x_j)+p_{*j})$$ such that
 $$z_*\in \tilde{B}^{n_1+t(x_1)}_{m_1}(f^{-t(x_1)}(x_1),\theta):=\cap_{j=-m_1}^{n_1+t(x_1)}f^{-j}B(f^{j-t(x_1)}(x_1),\theta q(f^{j-t(x_1)}(x_1))^{-2})$$
and for $2\leq i\leq k,$ $$f^{\sum_{j=1}^{i-1} (n_j+t(x_j)+p_{*j})+\sum_{j=2}^{i} m_j}(z_*)\in \tilde{B}^{n_i+t(x_i)}_{m_i}(f^{-t(x_i)}(x_i),\theta)$$$$=\cap_{j=-m_j}^{n_i+t(x_i)}f^{-j}B(f^{j-t(x_i)}(x_i),\theta q(f^{j-t(x_i)}(x_i))^{-2}).$$

Let $p_i=p_{*i}+t(x_{i+1}),\,\,p=p_*$ where $x_{k+1}=x_1,$ then $$\sum_{i=1}^{k} p_{i}\leq \sum_{i=1}^{k} K_*(\eta,\theta,f^{-t(x_i)}(x_i),m_i,n_i+t(x_i))+\sum_{i=1}^{k}t(x_{i+1})$$
$$=\sum_{i=1}^{k}K(\eta,\theta,x_i,m_i,n_i+t(x_i)).\,\,\,\,\,\,\,\,\,\,\,\,\,\,\,\,\,\,\,\,\,\,\,\,\,\,
\,\,\,\,\,\,\,\,\,\,\,\,\,$$
Let $z=f^{t(x_1)}(z_*)$, then $z$ is the needed periodic point. More precisely,
$$z=f^{t(x_1)}(z_*)\in f^{t(x_1)}\tilde{B}^{n_1+t(x_1)}_{m_1}(f^{-t(x_1)}(x_1),\theta)\,\,\,\,\,\,\,\,\,\,\,\,\,\,\,\,\,\,\,\,\,\,\,\,\,\,\,\,$$
$$=f^{t(x_1)}\cap_{j=-m_1}^{n_1+t(x_1)}f^{-j}B(f^{j-t(x_1)}(x_1),\theta q(f^{j-t(x_1)}(x_1))^{-2})\,\,\,$$$$=\cap_{j=-m_1-t(x_1)}^{n_1}f^{-j}B(f^j(x_1),\theta q(f^j(x_1))^{-2})\subseteq \tilde{B}^{n_1}_{m_1}(x_1,\theta)$$
and similarly for $2\leq i\leq k,$ we have $$f^{\sum_{j=1}^{i-1} (n_j+p_{j})+\sum_{j=2}^{i} m_j}(z)=f^{\sum_{j=1}^{i-1} (n_j+p_{j})+\sum_{j=2}^{i} m_j+t(x_1)}(z_*)\,\,\,\,\,\,$$
$$=\,f^{t(x_i)}\circ f^{\sum_{j=1}^{i-1} (n_j+t(x_j)+p_{*j})+\sum_{j=2}^{i} m_j}(z_*)\,\,\,\,\,\,\,\,\,\,\,\,\,\,\,\,
\,\,\,\,\,\,\,\,\,\,\,\,\,\,\,\,\,\,\,\,\,\,\,\,\,\,\,\,\,\,\,\,\,\,\,\,\,\,\,\,\,\,$$
$$\in\, f^{t(x_i)}\tilde{B}^{n_i+t(x_i)}_{m_i}(f^{-t(x_i)}(x_i),\theta)\,\,\,\,\,\,\,\,\,\,\,\,\,\,\,\,\,\,\,\,\,\,
\,\,\,\,\,\,\,\,\,\,\,\,\,\,\,\,\,\,\,\,\,\,\,\,\,\,\,\,\,
\,\,\,\,\,\,\,\,\,\,\,\,\,\,\,\,\,\,\,\,\,\,\,\,\,\,\,\,\,\,$$
$$=\,f^{t(x_i)}\cap_{j=-m_i}^{n_i+t(x_i)}f^{-j}B(f^{j-t(x_i)}(x_i),\theta q(f^{j-t(x_i)}(x_i))^{-2})\,\,\,\,\,\,\,\,\,\,\,\,\,\,\,\,\,\,\,\,\,\,$$
$$=\,\cap_{j=-m_i-t(x_i)}^{n_i}f^{-j}B(f^j(x_i),\theta q(f^j(x_i))^{-2})\subseteq \tilde{B}^{n_i}_{m_i}(x_i,\theta).\,\,\,\,\,\,\,\qed$$
\medskip

\begin{Rem}\label{Rem-Proof-GNS}From the above discussion, in fact one also has a precise description of $p_{*i}:$ $$ p_{*i}\leq K_*(\eta,\theta,x_{i},m,n)+ K_*(\eta,\theta,x_{i+1},m,n)$$  and then $$ p_i\leq K(\eta,\theta,x_{i},m,n)+ K(\eta,\theta,x_{i+1},m,n).$$ In other words, $p_{*i},p_i$ only depends on the point $x_{i}$ and next point $x_{i+1}.$
\end{Rem}

\section{Proof of Theorem \ref{Thm:HypMeaDomSpl}}

To prove Theorem \ref{Thm:HypMeaDomSpl} we need the exponentially shadowing lemma in \cite{Tian}($C^1$ Pesin theory). Before that we introduce some notions. Given $x\in M$ and $n\in\mathbb{N}$, let
   $$\{x,\,n\}:=
 \{f^j(x)\,|\,\,j=0,\,1,\,\cdots,\,n\}.$$ In other words, $\{x,\,n\}$
represents the orbit segment from $x$ to $f^n(x)$ with length $n$. For a sequence of points
$\{x_i\}_{i=-\infty}^{+\infty}$ in $M$
 and a sequence of positive integers
 $\{n_i\}_{i=-\infty}^{+\infty}$, we call $\{x_i,\,n_i\}_{i=-\infty}^{+\infty}$
 a $\delta$-pseudo-orbit,
 if $d(f^{n_i}(x_i),\,x_{i+1})<\delta$ for all $i$.
  Given $\varepsilon>0$ and $\tau>0,$ we call a point $x\,\in
M$  an (exponentially) $(\tau,\varepsilon)$-shadowing point for a pseudo-orbit
$\big{\{}x_i,\,n_i\big{\}}_{i=-\infty}^{+\infty},$  if
$$d\big{(}f^{c_i+j}(x),f^j(x_i)\big{)}<\tau\cdot e^{-\min\{j,\,n_i-j\}\varepsilon},$$ $\forall\,\,
j=0,\,1,\,2,\,\cdots,\,n_i$ and $\forall\,\, i \in \mathbb Z$, where
$c_i$
 is defined as
  \begin {equation} \label{eq:respective-time-squens}c_i=\begin{cases}
 0,&\text{for }i=0\\
 \sum_{j=0}^{i-1}n_j,\,&\text{for }i>0\\
 -\sum_{j=i}^{-1}n_j,\,&\text{for }i<0.
\end{cases}
 \end {equation}

\begin{Lem}\label{Lem:ExpShadowing} Let us assume same conditions as in Theorem \ref{Thm:HypMeaDomSpl}. Then for each $\sigma>0,$ there exist a compact set $\Lambda_\sigma\subseteq M$, $\varepsilon_\sigma>0$ and $T_\sigma \in \mathbb{N}$ such that $\mu(\Lambda_\sigma)>1-\sigma$ and following \emph{(Exponentially) Shadowing Lemma} holds. For
$\forall\,\,\tau>0,$ there exists $\delta=\delta(\sigma,\tau)>0$ such that if  a
$\delta$-pseudo-orbit $\{x_i,\,n_i\}_{i=-\infty}^{+\infty}$
satisfies $n_i\geq T_\sigma $ and $x_i,f^{n_i}(x_i)\in
\Lambda_\sigma$ for all $i$, then there exists an (exponentially)
$(\tau,\varepsilon_\sigma)$-shadowing point $x\in M $ for
$\{x_i,n_i\}_{i=-\infty}^{+\infty}$. If further
$\{x_i,n_i\}_{i=-\infty}^{+\infty}$ is periodic, i.e., there
 exists an integer $m>0$ such that $x_{i+m}=x_i$ and $n_{i+m}=n_i$ for all i,
 then the shadowing point $x$ can be chosen to be periodic.
\end{Lem}

\bigskip

{\bf Proof of Theorem \ref{Thm:HypMeaDomSpl}} Here we point out that the result of Lemma \ref{Lem:ExpShadowing} is  weaker than the statements of Katok's shadowing lemma since it holds only for the pseudo-orbit whose beginning and ending points in the same Pesin block. However,
 Lemma \ref{Lem:ExpShadowing} is  enough to prove Theorem \ref{Thm:HypMeaDomSpl}, since it also can deduce all propositions in Section 5. In other words, every ergodic measure of a homeomorphism with the property stated as  in Lemma \ref{Lem:ExpShadowing} has (generalized) non-uniform specification.
Here we only give a proof of  non-uniform specification(Definition \ref{Def:NS}).

Since the given hyperbolic measure $\mu$ is ergodic, the number $\varepsilon_\sigma$ in Lemma \ref{Lem:ExpShadowing} can be chosen independent on $\sigma$ from Remark 1.4 in \cite{Tian} and thus we can take a fixed number $\varepsilon$. In other word, for each $\sigma>0,$ there exist a compact set $\Lambda_\sigma\subseteq M$ and $T_\sigma \in \mathbb{N}$ such that $\mu(\Lambda_\sigma)>1-\sigma$ and  \emph{(Exponentially) Shadowing Lemma} holds as in Lemma \ref{Lem:ExpShadowing} for $\varepsilon_\sigma\equiv \varepsilon$.

 Set $\tilde \Lambda_\sigma=supp(
\mu|_{\Lambda_{\sigma}})$ and $\tilde
\Lambda=\cup_{\sigma>0}\tilde \Lambda_\sigma.$ Clearly,   $\tilde \Lambda$ is of
$\mu $-full measure.  Let $\Delta_{\sigma}\subseteq\tilde \Lambda_{\sigma}$ be the set of all points whose recurrence times are well defined for $\Gamma=\tilde \Lambda_{\sigma}$ and  satisfy
\begin {equation}\label{Eq:20Rec2}
\lim_{i,j\rightarrow+\infty}\frac{t_{i+1}(x)-t_{-j-1}(x)}{t_i(x)-t_{-j}(x)}=1.
\end {equation}
By Poincar\'{e} Recurrence Theorem and Proposition \ref{Prop:Rec2}, $\mu(\Delta_{\sigma})=\mu(\tilde \Lambda_{\sigma})$ and thus $\cup_{\sigma>0}\Delta_{\sigma}$ is a set with full measure. So we only need to prove that for every fixed $\Delta_{\sigma}$ with positive measure, all points in $\Delta_{\sigma}$ satisfy the conditions of non-uniform specification property.
Take and fix a point $x\in\Delta_{\sigma}$.

Let $0<\eta\leq {\varepsilon}/2$ and $q$ be an $\eta-$slowing varying positive function, $\theta>0$ and let $m,\,n$ be two positive integers. We may assume that $m,n\geq T_\sigma$(Otherwise, consider $m'=m+T_\sigma,\,n'=n+T_\sigma$).
Let $\tau=\theta q^{-2}(x)
>0$. Then for this $\tau$ there exists
$\delta=\delta(\tau,\sigma)>0$ satisfying Lemma \ref{Lem:ExpShadowing}.

Take and fix for $\tilde\Lambda_{\sigma}$
a finite cover $\alpha_\sigma=\{V_1,V_2,\cdots,V_{r_{\sigma}}\}$ by
nonempty open balls $V_i$ in $M$ such that diam$(U_i)
<\delta$ and $\mu(U_i)>0$ where $U_i=V_i\cap
\tilde\Lambda_{\sigma},\,i=1,\,2,\,\cdots,\,r_{\sigma}$. Since $\mu$ is $f-$ergodic, by Birkhoff ergodic theorem we have
\begin {equation}\label{Eq:21}
\lim_{n\rightarrow +\infty}\frac 1n
  \sum_{h=0}^{n-1}\mu(f^{-h}(U_i)\cap U_j)=\mu(U_i)\mu(U_j)>0.
\end{equation}
 Then
take
$$ X_{i,\,j}=\min\{h\in \mathbb{N}\,\,|\,\, h\geq
 \textbf{T}_\sigma, \,\,\mu(f^{-h}(U_i)\cap
  U_j)>0\}.
$$
By  (\ref{Eq:21}), $\textbf{T}_\sigma\leq X_{i,\,j}<+\infty$. Let
$$M_{\sigma}=
\max_{1\leq i, j\leq r_{\sigma}}X_{i,\,j}.$$
Note that $M_{\sigma}$ is dependent on $\sigma,\theta$ and the $\eta-$slowing varying positive function $q,$ but independent of $m,n.$ So \begin {equation}\label{Eq:22M_k} \frac{M_\sigma}m,\,\,\,\frac{M_\sigma}n\rightarrow 0
\end {equation} as $m,n\rightarrow \infty$ respectively.

Take positive integers $l_1$ and $l_2$  such that
\begin {equation}\label{Eq:22}
t_{-l_1}< -m \leq t_{-l_1+1}\,\,\,\,\text{and}\,\,\,\,t_{l_2}>n \geq t_{l_2-1}.
\end{equation}
Take positive integers $s_1$ and $s_2$  such that
 \begin {equation}\label{Eq:23}
 t_{-l_1-s_1}\leq (1+ \frac{2\eta} {\varepsilon}) t_{-l_1}<t_{-l_1-s_1+1}\,\,\,\,\text{and}\,\,\,\,t_{l_2+s_2}\geq (1+ \frac{2\eta} {\varepsilon}) t_{l_2}>t_{l_2+s_2-1}.
\end{equation}
Take $K:=K(\eta,\theta,x,m,n)=t_{l_2+s_2}-t_{-l_1-s_1}+M_\sigma-m-n.$ The calculation of $$\lim_{\eta\rightarrow 0}\limsup_{m,n\rightarrow+\infty}\frac{K(\eta,\theta,x,m,n)}{m+n}=0$$ is similar as in the
proof of Theorem \ref{Thm:NS} by using (\ref{Eq:22}), (\ref{Eq:22M_k}), (\ref{Eq:20Rec2}) and (\ref{Eq:23}). We omit the details.

Since $f^{t_{-l_1-s_1}}(x), f^{t_{l_2+s_2}}(x)\in\tilde\Lambda_{\sigma},$ we can take $U_i$ and $U_j$ such that
$$f^{t_{-l_1-s_1}}(x)\in U_i, f^{t_{l_2+s_2}}(x)\in U_j.$$ By  (\ref{Eq:21}), there exist $y\in U_j$ and $\textbf{T}_\sigma\leq N\leq M_\sigma$ such that $f^N(y)\in U_i.$

Note that $$f^{t_{-l_1-s_1}}(x),x, f^{t_{l_2+s_2}}(x), y, f^N(y)\in \tilde\Lambda_k\subseteq \Lambda_k,$$
$$-t_{-l_1-s_1}\geq m \geq T_\sigma,\,\,t_{l_2+s_2}\geq n \geq T_\sigma,\,\,N\geq T_\sigma$$ and $$d(f^{t_{l_2+s_2}}(x), y)<\delta,\,\,d(f^{t_{-l_1-s_1}}(x), f^N(y))<\delta.$$ So if we repeat the orbit segments of $$\{f^{t_{-l_1-s_1}}(x), -t_{-l_1-s_1}\},\,\,\,\{x,t_{l_2+s_2}\},\,\,\,\{y,N\}$$ infinite times, then we get a
periodic $\delta$ pseudo-orbit. Then by Lemma \ref{Lem:ExpShadowing} there exists a periodic point $z$ with period $p=t_{l_2+s_2}-t_{-l_1-s_1}+N$ such that
$$d\big{(}f^{j}(x),f^j(z)\big{)}<\tau\cdot e^{-\min\{j,\,t_{l_2+s_2}-j\}\varepsilon},$$ $\forall\,\,
j=0,\,1,\,2,\,\cdots,\,t_{l_2+s_2}$ and $$d\big{(}f^{j}(x),f^j(z)\big{)}<\tau\cdot e^{-\min\{-j,\,-t_{-l_1-s_1}+j\}\varepsilon},$$ $\forall\,\,
j=t_{-l_1-s_1},\,\,\cdots,\,-2,\,-1,\,0.$

Now we start to verify the conditions of non-uniform specification property. Clearly we have $p\leq m+n+K$ since $N\leq M_\sigma,$ and thus we only need to show $$z\in\tilde{B}^n_m(x,\theta):=\cap_{i=-m}^nf^{-i}B(f^k(x),\theta q(f^k(x))^{-2}).$$  Firstly, we consider $0\leq i\leq t_{l_2}-1$ and calculate $d(f^i(x),f^i(z)).$ More precisely,
\begin{eqnarray*}
& &d(f^i(x),f^i(z))\nonumber\\
&\leq& \max\{\tau e^{-i\varepsilon}, \tau  e^{-(t_{l_2+s_2}-i)\varepsilon}\}\nonumber\\
&\leq& \max\{\tau  e^{-i\varepsilon}, \tau  e^{-(t_{l_2+s_2}-t_{l_2})\varepsilon}\}\,\,\,\,\,\,\,(\text{using } i\leq t_{l_2})\nonumber\\
&\leq& \max\{\tau  e^{-i\varepsilon}, \tau  e^{-2t_{l_2}\eta}\}\,\,\,\,\,\,\,(\text{using } (\ref{Eq:23}))\nonumber\\
&\leq& \max\{\tau  e^{-2i\eta}, \tau  e^{-2i\eta}\}\,\,\,\,\,\,\,(\text{using } 2\eta<\varepsilon\text{ and }i\leq t_{l_2} )\nonumber\\
&=& \tau  e^{-2i\eta}\nonumber\\
&=& \theta q^{-2}(x)e^{-2i\eta}\,\,\,\,\,\,\,(\text{by the choice of } \tau )\nonumber\\
&\leq& \theta q(f^i(x))^{-2}\,\,\,\,\,\,\,\,(\text{using }q(f^i(x))\leq q(x)e^{i\eta}).\nonumber
\end{eqnarray*}
Secondly we can follow the similar method to show that for $t_{-l_1}+1\leq i\leq 0,$
\begin{eqnarray*}
d(f^i(x),f^i(z))\leq \theta q(f^i(x))^{-2}.\nonumber
\end{eqnarray*}
Notice that $t_{-l_1}<-m$ and $n< t_{l_2}$ and thus $z\in\tilde{B}^n_m(x,\theta)$.
\qed
\bigskip

\begin{Rem}\label{Rem-Further-NS} From the proofs of Theorem \ref{Thm:HypMeaDomSpl} and Theorem \ref{Thm:NS}, we point out that for  an invariant measure of a
homeomorphism $f:X\rightarrow X$   on a compact metric space, a sufficient condition of
non-uniform specification for $\mu$ is that: There exists $\varepsilon>0$ such that  for any $\sigma>0,$ there is a subset $\Gamma_\sigma$ with $\mu$ positive measure larger than $1-\sigma$ such that for  any $\tau>0$  and any integer $m,\,n,$ if $f^{-m}(x),x,f^n(x)\in\Gamma_\sigma$ there exists $L:=L(\sigma,\tau,x,m,n)$ such that:\\
(i) there exists a periodic point $z$ with period $p\leq m+n+L$ such that
$$d\big{(}f^{j}(x),f^j(z)\big{)}<\tau\cdot e^{-\min\{j,\,n-j\}\varepsilon},$$ $\forall\,\,
j=0,\,1,\,2,\,\cdots,\,n$ and $$d\big{(}f^{j}(x),f^j(z)\big{)}<\tau\cdot e^{-\min\{-j,\,m+j\}\varepsilon},$$ $\forall\,\,
j=-m,\,\,\cdots,\,-2,\,-1,\,0;$
\\
(ii) the dependence of $L$ on $m,\,n$ satisfies
$$\limsup_{m,n\rightarrow+\infty}\frac{L(\sigma,\tau,x,m,n)}{m+n}=0.$$\\
In particular, exponentially shadowing is such a sufficient condition. Note that from the proofs of Theorem \ref{Thm:HypMeaDomSpl} and Theorem \ref{Thm:NS},  $L$ is the number $M_k$ or $M_\sigma$  independent of $m,n$. We can also give a similar sufficient condition with several orbit segments for generalized non-uniform specification(we omit the details). There are also some other related papers\cite{WS,Dai,Kal,BP,Tian} that exponential shadowing plays important roles. Exponential closing (which is the particular exponential shadowing for one pseudo orbit segment) has played crucial roles in \cite{WS,Dai,Kal}  to prove that Lyapunov exponents of ergodic measures can be approximated by ones of periodic measures and it has  also been used in \cite{BP,Tian,Kal} to calculate   H$\ddot{\text{o}}$lder functions or H$\ddot{\text{o}}$lder cocycles to get some convergence properties for proving corresponding Liv$\check{\text{s}}$ic Theorem.
\end{Rem}

{\bf Acknowledgement.} The second author thanks very much
to Professor Marcus Bronzi  and Fernando Pereira Micena {\it et al} for their help in the period of  UFAL.


\section*{ References.}
\begin{enumerate}


\itemsep -2pt

\small


\bibitem{Bl83} A.M. Blokh, Decomposition of dynamical systems on an interval, Russ. Math. Surv. 38, 133--134 (1983).

\bibitem{Bow}R. Bowen. Periodic orbits for hyperbolic flows, Amer. J. Math., 94 (1972), 1-30.

\bibitem{Bowen2}R. Bowen, Equilibrium states and the ergodic theory of Anosov
diffeomorphisms, Springer Lecture Notes in Math. 470 (1975).

\bibitem{BP} Luis Barreira and Yakov B. Pesin,  Nonuniform hyperbolicity, Cambridge Univ. Press, Cambridge (2007).

\bibitem{Bu97} J. Buzzi, Specification on the Interval, Transactions of the American Mathematical Society,
Volume 349, Number 7,  Pages 2737--2754 (1997).

\bibitem{Dai} X. Dai,  Exponential closing property and approximation of Lyapunov exponents of linear cocycles, Forum. Math., to appear.


\bibitem{Hir} M. Hirayama, Periodic probability measures are dense in
the set of invariant measures,  Dist.  Cont. Dyn. Sys. 9 ( 2003),
1185-1192.

\bibitem{Kal} B. Kalinin,  Liv$\check{s}$ic Theorem for matrix cocycles, to appear in Annals of Mathematics.


\bibitem{Katok1}A. Katok and B. Hasselblatt, Introduction to the modern theory of
dynamical systems, Encyclopedia of Mathematics and its Applications
54, Cambridge Univ. Press, Cambridge (1995).

\bibitem{Katok2}A. Katok, Lyapunov
exponents, entropy and periodic orbits for diffeomorphisms, Inst.
Hautes Etudes Sci. Publ. Math. 51 (1980), 137-173.

\bibitem{LLS} C. Liang, G. Liu, W. Sun,  Approximation properties
on invariant measure and Oseledec splitting in non-uniformly
hyperbolic systems,  Trans. Amer. Math. Soci.  361 (2009) 1543-1579.


\bibitem{LST} C. Liang,  W. Sun, X. Tian, Ergodic properties of invariant measures for $C^{1+\alpha}$ non-uniformly hyperbolic systems, arXiv:1011.5300.

\bibitem{LST1} G. Liao,  W. Sun, X. Tian,  Metric entropy and the number of periodic points,  Nonlinearity  23 (2010) 1547-1558.


\bibitem{Mar} B. Marcus, A note on periodic points for ergodic toral automorphisms, Monatsh. Math., 89,
121-129 (1980).


\bibitem{Oli} K. Oliveira, Every expanding measure has the non-uniform specification property, arXiv: 1007.1449.


\bibitem{OV08} K. Oliveira, M. Viana, Thermodynamical formalism for robust classes of potentials and non-uniformly hyperbolic maps, Ergodic Theory \& Dynamical Systems (Print), v. 28, p. 501-533, 2008.


\bibitem{Os}V. I. Oseledec, A multiplicative ergodic theorem, Trans. Mosc. Math. Soc, 19 (1968), 197-231.

\bibitem{Pollicott}M. Pollicott, Lectures on ergodic theory and
Pesin theory on compact manifolds, Cambridge Univ. Press, Cambridge
(1993).

\bibitem{PS} C. Pfister and W. Sullivan, On the topological entropy of saturated sets, Ergodic Theory and Dynamical Systems, 27 (2007), 929-956.


\bibitem{STV} B. Saussol, S. Troubetzkoy and S. Vaienti, Recurrence and Lyapunov exponents, Mosc. Math. J. , Volume 3, Number 1 (2003) 189-203.

\bibitem{Sigmund}K. Sigmund, Generic properties of invariant measures for axiom A
diffeomorphisms, Invention Math. 11(1970), 99-109.

\bibitem{ST} W. Sun, X. Tian,  Pesin set, closing lemma and shadowing lemma in $C^1$
non-uniformly hyperbolic systems  with limit domination, arXiv:1004.0486.


\bibitem{T} D. Thompson, A variational principle for topological pressure for certain non-compact sets, Journal of the London Mathematical Society, Volume 80, part 3 (2009), 585-602.

\bibitem{Tian} X. Tian, Hyperbolic measures with dominated   splitting, arXiv:1011.6011v2.



\bibitem{Var} P. Varandas, Non-uniform specification and large deviations for weak Gibbs measures, arXiv:0906.3350v2.


\bibitem{WS}Z. Wang and W. Sun, Lyapunov exponents of hyperbolic measures and hyperbolic period
orbits, Trans. Amer. math. Soc., 362  (2010), 4267-4282.

\bibitem{Y} K. Yamamoto, On the weaker forms of the specification property and their applications, Proceedings of the American Mathematical Society, V. 137, 11 (2009)  3807-3814.




\end{enumerate}
\end{document}